\documentclass[12pt]{amsart}                 %
 \theoremstyle{plain}    
 \newtheorem{thm}{Theorem}[section]
 \numberwithin{equation}{section} 
 \numberwithin{figure}{section} 
 \theoremstyle{plain}    
 \newtheorem*{thm*}{Theorem} 
 \theoremstyle{plain}    
 \newtheorem*{cor*}{Corollary}
 \theoremstyle{plain}    
 \newtheorem{lem}[thm]{Lemma} 
 \theoremstyle{plain}
 \newtheorem*{lem*}{Lemma}
 \theoremstyle{plain}    
 \newtheorem{prop}[thm]{Proposition} 
 \theoremstyle{plain}
 \newtheorem*{prop*}{Proposition}
 \theoremstyle{remark}
 \newtheorem{rem}[thm]{Remark}
 \theoremstyle{remark}    
 \newtheorem*{acknowledgement*}{Acknowledgement} 

\usepackage{amssymb}
\usepackage{eepic}

\def\id{\operatorname{id}}
\def\diag{\operatorname{diag}}
\begin{document}

\title{Irreducible subfactors of \protect\( L(\mathbb {F}_{\infty })\protect \) of
index \protect\( \lambda >4\protect \).}

\author{Dimitri Shlyakhtenko and Yoshimichi Ueda}

\address{Department of Mathematics, UCLA, Los Angeles, CA 90095, USA.}

\email{shlyakht@math.ucla.edu}

\address{Department of Mathematics, Graduate School of Science, Hiroshima University,
Higashi-Hiroshima, 739-8526, Japan.}

\email{ueda@math.sci.hiroshima-u.ac.jp }

\begin{abstract}
By utilizing an irreducible inclusion of type III\( _{q^{2}} \) factors coming
from a free-product type action of the quantum group \( SU_{q}(2) \), we show
that the free group factor \( L(\mathbb {F}_{\infty }) \) possesses irreducible
subfactors of arbitrary index \( >4 \). Combined with earlier results of Radulescu,
this shows that \( L(\mathbb {F}_{\infty }) \) has irreducible subfactors with
any index value in \( \{4\cos ^{2}(\pi /n):n\geq 3\}\cup [4,+\infty ) \).
\end{abstract}
\maketitle

\section{Introduction.}

Let \( M \) be a type II\( _{1} \) factor, and let \( N\subset M \) be a
finite index subfactor. Recall that \( N\subset M \) is called \emph{irreducible,}
if the relative commutant \( N'\cap M \) \emph{}is trivial. 

In his fundamental paper \cite{jones:index}, V.F.R.~Jones introduced the notion
of index \( [M:N] \) of a subfactor \( N\subset M \). Related to the index,
he introduced two invariants of a factor \( M \):
\begin{eqnarray*}
\mathcal{I}(M) & = & \{\lambda \in \mathbb {R}_{+}:\exists N\subset M\textrm{ subfactor},\, [M:N]=\lambda \},\\
\mathcal{C}(M) & = & \{\lambda \in \mathbb {R}_{+}:\exists N\subset M\textrm{ irreducible}\, \textrm{subfactor},\, 
[M:N]=\lambda \},
\end{eqnarray*}
measuring the possible index values of subfactors of \( M \). 

In the same paper, Jones showed that for any factor \( M \), the index of a
subfactor cannot be arbitrary; in fact, \( \mathcal{C}(M)\subset \mathcal{I}(M)\subset \mathcal{I}=\{4\cos ^{2}(\pi /n):n\geq 
3\}\cup [4,+\infty ] \).
Moreover, he produced subfactors of the hyperfinite II\( _{1} \) factor \( R \)
of arbitrary allowable index, thus showing that \( \mathcal{I}(R)=\mathcal{I} \).
However, the subfactors of \( R \) with index \( >4 \) that he constructed
are not irreducible; to this day, it is still not known whether the equality
\( \mathcal{C}(R)=\mathcal{I} \) holds; in other words, whether every index
value can be realized by an irreducible subfactor of the hyperfinite factor.
The one general positive result is due to by Popa \cite{popa:markov} (see also
\cite{popa:universalconstructions} and \cite{popa:standardlattice}), stating
that given a number \( \lambda \in \mathcal{I} \), there is a (non-hyperfinite)
II\( _{1} \) factor \( M \), for which \( \lambda \in \mathcal{C}(M) \). 

Thanks to the theory of free probability, pioneered by Voiculescu in the early
80s (see e.g. \cite{DVV:book}), there has been much progress at understanding
free group factors \( L(\mathbb {F}_{n}) \) (Voiculescu's introduction to \cite{DVV:book}
states that these are the ``best'' of the ``bad'' non-hyperfinite factors).
Following this philosophy, Radulescu undertook a study of subfactors of free
group factors. Using a particular version of Popa's construction \cite{radulescu:subfact}
and by employing free probability and random matrix techniques, he managed to
show that for an interpolated free group factor \( L(\mathbb {F}_{t}) \), \( \mathcal{C}(L(\mathbb {F}_{t}))\cap 
[0,4]=\mathcal{I}(L(\mathbb {F}_{t}))\cap [0,4]=\mathcal{I}\cap [0,4] \).
Coupled with his earlier result \cite{radulescu2} on the fundamental group
of \( L(\mathbb {F}_{\infty }) \) this gave that \( L(\mathbb {F}_{\infty }) \)
has (possibly non-irreducible) subfactors of all indices, i.e., \( \mathcal{I}(L(\mathbb {F}_{\infty }))=\mathcal{I} \).
However, this still left the question about whether subfactors of \( L(\mathbb {F}_{\infty }) \) of index \( >4 \) can be chosen to be irreducible.

The main theorem of this paper is

\begin{thm*}
For each \( \lambda >4 \), there exists an irreducible subfactor of \( L(\mathbb {F}_{\infty }) \)
of index \( \lambda  \).
\end{thm*}
Combined with earlier results of F. Radulescu, this gives:

\begin{cor*}
\( \mathcal{C}(L(\mathbb {F}_{\infty }))=\mathcal{I}(L(\mathbb {F}_{\infty }))=\mathcal{I} \).
\end{cor*}
In fact, more is true: given an allowable index value \( \lambda \in \mathcal{I} \),
there exists an irreducible subfactor \( L(\mathbb {F}_{\infty })\cong N\subset M\cong L(\mathbb {F}_{\infty }) \)
of index \( \lambda  \).

We would like to point out that it is still not known whether every Popa's \( \lambda  \)-lattice
\cite{popa:standardlattice} (i.e., an abstract system of higher relative commutant
of a subfactor) can be realized by a subfactor of some fixed universal II\( _{1} \)
factor \( M \). The factor \( M=L(\mathbb {F}_{\infty }) \) seems to be a
likely candidate, although we don't know how to prove it; note, however, that
Radulescu gave an affirmative answer to this question in the case that the \( \lambda  \)-lattice
is finite-depth. Also, by the results of this paper, all examples of irreducible
subfactors whose \( \lambda  \)-lattices come from the representation theory
of \( SU_{q}(2) \) (see \cite{banica:compactQuantumGroupsSubfactors}, \cite{xu:latticesFromQuantumGroups},
\cite{wenzl:tensorCatergoriesQuantumGroups}) can be realized as \( N\subset M \)
with \( N\cong M\cong L(\mathbb {F}_{\infty }) \). 

We don't know whether our result holds for \( L(\mathbb {F}_{n}) \), \( n<\infty  \).

The strategy to find irreducible subfactors of \( L(\mathbb {F}_{\infty }) \)
with index \( \lambda >4 \) comes from the work of the second-named author
\cite{ueda:fixed-point} on quantum \( SU_{q}(2) \) actions on free products
of von Neumann algebras, and analysis of resulting Wassermann-type inclusions
of type III factors. As was pointed out in \cite[Question 2]{ueda:fixed-point},
to prove the main theorem of the present paper, it is sufficient to identify
a certain crossed product of a certain von Neumann algebra by \( SU_{q}(2) \),
which we do in Section \ref{sec:amalgfreeprod} by identifying this algebra
as a free Araki-Woods factor \cite{shlyakht:quasifree:big}. For the convenience
of the reader, we summarize the main necessary facts from \cite{ueda:fixed-point}
in Section \ref{sec:quantumactions}.

\begin{acknowledgement*}
The authors would like to express their sincere gratitude to MSRI and the organizers
of the year-long program ``Operator Algebras (2000--2001)\char`\"{} held at
MSRI for inviting them to the program, where this joint work (although conceived
some time ago) was carried out and completed. 

D.S. would like to acknowledge support by a National Science Foundation Postdoctoral
Fellowship. Y.U. would like to acknowledge support by the Japanese Ministry
of Education, Science, Sports and Culture.
\end{acknowledgement*}

\section{Amalgamated free products of free group factors with \protect\( B(H)\protect \).\label{sec:amalgfreeprod}}

In this section, we prove some technical results about amalgamated free products
of free group factors with \( B(H) \), identifying them as free Araki-Woods
factors. The treatment given in this section is more general than strictly speaking 
necessary for our result, but allows for more succinct and abstract proofs.
We encourage the interested reader to look through Appendix I, which can serve as 
an illustration of the proof given here. 

We recall some notation from \cite{shlyakht:semicirc}. Let \( M \) be a von
Neumann algebra and \( \eta :M\to M\otimes B(\ell ^{2}) \) be a normal completely
positive map, given as a matrix \( \eta (m)=(\eta _{ij}(m))_{ij} \). Let \( K \)
be the \( M \)-Hilbert bimodule spanned by vectors \( \xi _{i} \), and satisfying
\[
\langle \xi _{i},m\xi _{j}\rangle _{M}=\eta _{ij}(m),\quad \forall m\in M\]
(see \cite[Lemma 2.2]{shlyakht:semicirc}). 

Let \( \mathcal{T} \) be the universal \( C^{*} \)-algebra generated by \( M \)
and operators \( L(\zeta ):\zeta \in K \) satisfying
\[
L(\zeta )^{*}mL(\zeta ')=\langle \zeta ,m\zeta '\rangle ,\quad \forall m\in M,\, \zeta ,\zeta '\in K.\]
This is the Toeplitz extension of the Cuntz-Pimsner \( C^{*} \)-algebra associated
to the \( M \)-bimodule \( H \) (see \cite{pimsner}). The operators \( L(\zeta ) \)
are called \( \eta  \)-creation operators over \( M \).

Set \( L_{i}=L(\xi _{i}) \). Let \( E_{M} \) denote the canonical conditional
expectation from \( \mathcal{T} \) onto \( M \) (see \cite{pimsner}, \cite{shlyakht:semicirc}).
Choose a faithful normal state \( \phi _{0} \) on \( M \) and set \( \phi =\phi _{0}\circ E_{M} \).
In the GNS representation associated to \( \phi  \), consider the von Neumann
algebra
\[
\Phi (M,\eta )=W^{*}(M,L_{i}+L_{i}^{*}:i=1,2,\ldots ).\]
It can be shown (\cite{shlyakht:semicirc}) that \( \Phi (M,\eta ) \) depends
only on \( M \) and \( (\eta _{ij}) \). 

Let \( K^{\eta }\subset K \) be the \( \mathbb {R} \)-linear subspace of the
Hilbert \( M \)-bimodule \( K \), given as the norm closure of the set of
vectors
\[
\sum _{\textrm{finite}}m_{j}\xi _{j}n^{*}_{j}+n_{j}\xi _{j}m_{j}^{*},\quad m_{j},n_{j}\in M.\]
Then \cite[Proposition 2.19]{shlyakht:semicirc} \( \Phi (M,\eta ) \) depends
only on the inclusion \( K^{\eta }\subset K \). In fact
\begin{equation}
\label{eqn:PhiandKeta}
\Phi (M,\eta )=W^{*}(M,L(\zeta )+L(\zeta )^{*}:\zeta \in K^{\eta }\subset K)
\end{equation}
(in the GNS representation associated to \( \phi  \)).

In the particular case that \( M=\mathbb {C} \), the matrix \( \eta _{ij}(1) \)
defines a scalar-valued inner product on \( K \), and \( K^{\eta }\subset K \)
ends up a particular real subspace. In this case, under the assumption that
the state \( E_{M}=E_{\mathbb {C}} \) is faithful, \( \Phi (M,\eta ) \) is
nothing other than the free Araki-Woods factor \( \Gamma (K^{\eta }\subset K)'' \)
in the notation of Remark 2.6 of \cite{shlyakht:quasifree:big} (see also \cite[Example 3.5]{shlyakht:semicirc}). 

A very particular example of such a matrix \( \mu  \) is
\[
\mu _{\lambda }=\left( \begin{array}{cc}
1 & -i\frac{\lambda -1}{\lambda +1}\\
i\frac{\lambda -1}{\lambda +1} & 1
\end{array}\right) ,\quad 0<\lambda <1;\]
in this case \( \Phi (\mathbb {C},\mu _{\lambda }) \) is of type III\( _{\lambda } \)
and is the unique free Araki-Woods factor \( T_{\lambda } \) of this type. 

\begin{thm}
\label{thrm:genresultA-valued}Let \( H \) be a finite or infinite dimensional
Hilbert space. Let \( (\eta _{ij}):B(H)\to B(H)\otimes B(\ell ^{2}) \) be any
normal completely positive map. Assume that \( E_{B(H)}:\Phi (B(H),\eta )\to B(H) \)
is faithful. Then the von Neumann algebra \( \Phi (B(H),\eta ) \) is stably
isomorphic to the von Neumann algebra \( \Phi (\mathbb {C},\mu ) \) for some
(completely) positive map \( \mu :\mathbb {C}\to B(\ell ^{2}) \); moreover,
\( E_{\mathbb {C}}:\Phi (\mathbb {C},\mu )\to \mathbb {C} \) is faithful.
\end{thm}
\begin{proof}
By equation (\ref{eqn:PhiandKeta}),
\[
N=\Phi (B(H),\eta )=W^{*}(B(H),L(\zeta )+L(\zeta )^{*}:\zeta \in K^{\eta }\subset K),\]
where \( K \) and \( K^{\eta } \) are as described above. Let \( e_{ij} \),
\( 0\leq i,j<\infty  \) be a family of matrix units generating \( B(H) \);
thus
\[
e_{ij}e_{j'i'}=\delta _{jj'}e_{ii'},\quad e_{ij}^{*}=e_{ji}.\]
Consider the algebra
\[
P=e_{00}Ne_{00}.\]
It is clear that \( P\otimes B(\ell ^{2})\cong N\otimes B(\ell ^{2}) \); thus
it is sufficient to prove that \( P \) is isomorphic to \( \Phi (\mathbb {C},\mu ) \)
for some \( \mu  \). Note that \( P \) is generated by the family
\[
e_{0i}(L(\zeta )+L(\zeta )^{*})e_{j0},\quad 0\leq i\leq j,\, \zeta \in K^{\eta };\]
hence
\[
P=W^{*}(L(e_{0i}\zeta e_{j0})+L^{*}(e_{0j}\zeta e_{i0}):\zeta \in K^{\eta }).\]
Consider the linear space
\begin{eqnarray*}
V_{\mathbb {R}}=\textrm{span}_{\mathbb {R}}(\{e_{0i}\zeta e_{i0}:\zeta \in K^{\eta },i\geq 0\}\cup \{e_{0i}\zeta 
e_{j0}+e_{0j}\zeta e_{i0}:\zeta \in K^{\eta },0\leq i<j\} &  & \\
\cup \{\sqrt{-1}(e_{0i}\zeta e_{j0}-e_{0j}\zeta e_{i0}):\zeta \in K^{\eta },0\leq i<j\})\subset K^{\eta }, &  & 
\end{eqnarray*}
as a subspace of the complex linear space
\[
V=\textrm{span}\{e_{0i}\zeta e_{j0}:\zeta \in K\}\subset K.\]
Since for \( i<j \),
\begin{eqnarray*}
(L(e_{0i}\zeta e_{j0})+L^{*}(e_{0j}\zeta e_{i0})) & = & \frac{1}{2}L(e_{0i}\zeta e_{j0}+e_{0j}\zeta e_{i0})+L(e_{0i}\zeta 
e_{j0}+e_{0j}\zeta e_{i0})^{*}\\
 &  & -\frac{1}{2}\sqrt{-1}(L(\sqrt{-1}(e_{0i}\zeta e_{j0}-e_{0j}\zeta e_{i0}))\\
 &  & \quad +L(\sqrt{-1}(e_{0i}\zeta e_{j0}-e_{0j}\zeta e_{i0}))^{*}),
\end{eqnarray*}
it follows that
\[
P=W^{*}(L(\zeta )+L(\zeta )^{*}:\zeta \in V_{\mathbb {R}}).\]
Since for any \( \zeta ,\zeta '\in V \), \( e_{00}\zeta e_{00}=\zeta  \) and
\( e_{00}\zeta 'e_{00}=\zeta ' \), we get that
\[
\langle \zeta ,\zeta '\rangle _{B(H)}=\langle e_{00}\zeta e_{00},e_{00}\zeta 'e_{00}\rangle =e_{00}\langle \zeta ,\zeta 
'\rangle e_{00}\in \mathbb{C} e_{00}.
\]
It follows that the restriction of the inner product on \( K \) to \( V \)
is scalar-valued; we denote the restriction by \( \langle \cdot ,\cdot \rangle _{V} \).
Hence we get an inclusion of a real Hilbert space into a complex Hilbert space
\[
V_{\mathbb {R}}\subset V.\]
In particular, for any \( \zeta ,\zeta '\in V \), we have
\begin{equation}
\label{eqn:cuntzrelations}
L(\zeta )^{*}L(\zeta ')=\langle \zeta ,\zeta '\rangle _{V}.
\end{equation}

Denote by \( \theta  \) the state \( P\ni m\mapsto e_{00}E_{B(H)}(m)e_{00} \).
Since \( E_{B(H)} \) is faithful by assumption, and \( \theta(e_{00}me_{00})=E_{B(H)}(e_{00}me_{00}) \),
we get that \( \theta  \) is a faithful state. Furthermore,
\begin{equation}
\label{eqn:thetaonLs}
\theta (L(\zeta _{1})\cdots L(\zeta _{k})L(\zeta _{k+1})^{*}\cdots L(\zeta _{r})^{*})=0
\end{equation}
if \( r\neq 0 \); note that this, together with the relations (\ref{eqn:cuntzrelations})
determines \( \theta  \) on \( P \).

Choose now a basis \( \zeta _{i} \) for \( V_{\mathbb {R}} \) as a real Hilbert
space, and let
\[
\mu _{ij}=\langle \zeta _{i},\zeta _{j}\rangle _{V}.\]
Then \( \Phi (\mathbb {C},\mu ) \) is generated in the GNS representation associated
to \( E_{\mathbb {C}} \) by \( l_{i}+l_{i}^{*} \), subject to the relations
\[
l_{i}^{*}l_{j}=\mu _{ij}=\langle \zeta _{i},\zeta _{j}\rangle _{V}.\]
Moreover, \( \psi =E_{\mathbb {C}} \) is determined by
\[
\psi (l_{i_{1}}l_{i_{2}}\cdots l_{i_{k}}l_{i_{k+1}}^{*}\cdots l_{i_{r}}^{*})=0,\quad r\neq 0.\]
Comparing these with (\ref{eqn:cuntzrelations}) and (\ref{eqn:thetaonLs}),
and using the fact that \( \theta  \) is faithful, we obtain that the map \( L(\zeta _{i})+L(\zeta _{i})^{*}\mapsto l_{i}+l_{i}
^{*} \)
extends to an isomorphism \( P\cong W^{*}(l_{i}+l_{i}^{*})=\Phi (\mathbb {C},\mu ) \),
which conjugates \( \theta  \) and \( \psi =E_{\mathbb {C}} \). It also follows
that \( E_{\mathbb {C}} \) is faithful.
\end{proof}
We can now deduce the main technical result needed for further computations;
for convenience, we write \( L(\mathbb {F}_{1})=L(\mathbb {Z}) \).

\begin{thm}
\label{thrm:freeproductBofH}Let \( H \) be a separable Hilbert space (finite
or infinite-dimensional). Let \( B\subset B(H) \) be a von Neumann subalgebra,
and assume that \( E:B(H)\to B \) is a normal faithful conditional expectation.
For \( m=1,2,\ldots  \) or \( +\infty  \), consider the reduced amalgamated
free product
\[
N=(L(\mathbb {F}_{m})\otimes B,\tau \otimes \id )*_{B}(B(H),E).\]
Then \( N \) is stably isomorphic to a free Araki-Woods factor \( \Gamma (\mathcal{H}_{\mathbb {R}},U_{t})'' \).

In particular, if \( N \) is of type III\( _{\lambda } \), \( 0<\lambda <1 \),
then \( N \) is isomorphic to the unique type III\( _{\lambda } \) free Araki-Woods
factor \( T_{\lambda } \).
\end{thm}
\begin{proof}
Let \( \eta :B(H)\to B(H)\otimes B(\ell ^{2}) \) be given by
\[
\eta (T)=\diag (E(T),E(T),\dots ,E(T),0,0\dots )\]
(\( m \) copies of \( E(T) \)). Since \( \eta_{ij} = \delta_{ij}E \), the
canonical conditional expectation \( E_{B(H)} : \Phi(B(H),\eta)
\rightarrow B(H) \) is known to be faithful (see e.g. Appendix II). 
By \cite[Example 3.3(c)]{shlyakht:semicirc},
\[
\Phi (B(H),\eta )\cong (B(H),E)*_{B}(\Phi (B,\eta |_{B}),E_{B}).\]
Moreover, by \cite[Example 3.3(b)]{shlyakht:semicirc},
\[
(\Phi (B,\eta |_{B}),E_{B})\cong (\Phi (\mathbb {C},\eta |_{\mathbb {C}1_{B}})\otimes B,E_{\mathbb {C}}\otimes \id )\cong 
(L(\mathbb {F}_{m})\otimes B,\tau \otimes \id ).\]
Hence
\[
\Phi (B(H),\eta )\cong (L(\mathbb {F}_{m})\otimes B,\tau \otimes \id
)*_{B}(B(H),E) = N.\] Applying Theorem \ref{thrm:genresultA-valued}, we
get that
\[
N\otimes B(\ell ^{2})\cong \Phi (\mathbb {C},\mu )\otimes B(\ell ^{2})\]
for some \( \mu  \).  It follows that \( \Phi (\mathbb {C},\mu )\cong \Gamma (\mathcal{H}_{\mathbb {R}},U_{t})'' \)
for some real Hilbert space \( \mathcal{H}_{\mathbb {R}} \) and one-parameter
group \( U_{t} \) (in fact, \( \mathcal{H}_{\mathbb {R}} \) is the real subspace
\( V_{\mathbb {R}} \) constructed in the proof of Theorem \ref{thrm:genresultA-valued}).

If \( N \) is of type III\( _{\lambda } \), \( 0<\lambda <1 \), it follows from
the uniqueness of the type III\( _{\lambda } \) free Araki-Woods factor \( T_{\lambda } \)
\cite{shlyakht:quasifree:big} that \( N\cong T_{\lambda } \).
\end{proof}
The following theorem can be deduced from the results of Radulescu \cite{radulescu1}
and Dykema \cite{dykema:typeiii}, as was done in \cite{shlyakht:quasifree:big}.
We have since found a somewhat shorter argument, which we give below:

\begin{thm}
\label{thrm:coreofTlambda}Let \( T_{\lambda } \) be the (unique) free Araki-Woods
factor of type III\( _{\lambda } \). Let \( \phi  \) be any normal faithful
state, so that the modular group \( \sigma ^{\phi }_{t} \) is periodic of period
\( - 2\pi /\log \lambda  \). Then the centralizer \( T_{\lambda }^{\phi } \)
is isomorphic to \( L(\mathbb {F}_{\infty }) \). Also, the core of \( T_{\lambda } \)
is isomorphic to \( L(\mathbb {F}_{\infty })\otimes B(H) \).
\end{thm}
\begin{proof}
Clearly, only the statement about the centralizer needs to be proved, since
the core is isomorphic to the centralizer, tensor \( B(H) \), as soon as the
centralizer is a factor (see \cite{connes}). Also, we can restrict ourselves
to dealing with a particular choice of the state \( \phi  \), satisfying the
hypothesis of this theorem; indeed, by \cite{connes}, all centralizers of such
states are stably isomorphic. Since by Radulescu's results \cite{radulescu2},
the fundamental group of \( L(\mathbb {F}_{\infty }) \) is all of \( (0,+\infty ) \),
a II\( _{1} \) factor \( M \) is stably isomorphic to \( L(\mathbb {F}_{\infty }) \)
iff it is actually isomorphic to \( L(\mathbb {F}_{\infty }) \).

Recall \cite{shlyakht:quasifree:big} that \( T_{\lambda } \) can be  described
as the free product
\[
L^{\infty }[0,1]*(B(H),\theta )\]
where \( \theta  \) is a normal faithful  state on \( B(H) \) given by
\[
\theta (T)=\text {Tr}(DT),\]
where
\[
D=\text {diag}(\frac{1}{1-{\lambda }},\frac{{\lambda }}{1-{\lambda }},\cdots ).
\]

Consider now the completely-positive  map \( \eta  \)  given by \( \eta (T)=\theta (T) \),
and consider \( \Phi (B(H),\eta )=W^{*}(B(H),L+L^{*}) \), where \( L^{*}TL=\eta (T)=\theta (T) \)
for all \( T\in B(H) \).  By \cite[Example 3.3(a)]{shlyakht:semicirc}, \( \Phi (B(H),\eta )\cong (B(H),\theta )*(L^{\infty }[0,
1],\textrm{Lebesgue measure})\cong (T_{\lambda },\phi _{\lambda }) \). 

Note that \( B(H) \) is generated, as a von Neumann algebra, by a partial isometry
\( v \), so that \( v^{*}v=1 \) and \( vv^{*}=\text {diag}(0,1,1,\ldots ) \).
Hence \( T_{\lambda } \) is generated by \( L+L^{*} \) and \( v \).  The
modular group of \( \phi _{\lambda } \) acts by fixing \( L+L^{*} \) and sending
\( v \) to \( {\lambda }^{it}v \).  The centralizer of the free product state
is then isomorphic to the von Neumann algebra generated by all elements of
the form
\[
v^{k}(L+L^{*})(v^{*})^{k},\quad (v^{*})^{l}(L+L^{*})v^{l},\quad k,l=0,1,\dots \]
as well as the projections \( p_{k}=v^{k}(v^{*})^{k} \). 

Let for \( k\geq 0 \)
\[
L_{k}=v^{k}L(v^{*})^{k}\]
and for \( -l<0 \),
\[
L_{-l}=(v^{*})^{l}Lv^{l}.\]
Denote by \( \Psi  \) the faithful normal conditional expectation
\[
\Psi (m)=\int _{0}^{-2\pi/\log\lambda }\sigma _{t}^{\phi _{\lambda }}(m)dt,\quad \Psi :T_{\lambda }\to T_{\lambda }^{\phi 
_{\lambda }}.\]
It follows that any \( m\in T_{\lambda } \) can be written as the sum (convergent
in \( L^{2}(T_{\lambda },\phi _{\lambda }) \))
\[
m=\sum _{k\geq 0}v^{k}m_{k}+\sum _{k<0}m_{k}(v^{*})^{k},\quad m_{k}=\Psi ((v^{*})^{-k}m),k\geq 0,m_{k}=\Psi (mv^{-k}),k<0,\]
with \( m_{k}\in T_{\lambda }^{\psi _{\lambda }} \). It follows that the centralizer
is linearly spanned by words of the form
\[
v^{k_{1}}(v^{*})^{l_{1}}Xv^{k_{2}}(v^{*})^{l_{2}}X\cdots v^{k_{p}}(v^{*})^{l_{p}},\quad \sum k_{j}-\sum l_{j}=0.\]
where \( X=L+L^{*} \). Using the notation
\[
X_{k}=L_{k}+L_{k}^{*},\quad k\in \mathbb {Z}\]
and the relation \( v^{*}v=1 \), \( v^{k}(v^{*})^{k}=p_{k}\in D \), we get
that the centralizer is generated by
\[
D\textrm{ and }\{X_{k}:k\in \mathbb {Z}\}.\]

Note that for all diagonal elements \( d\in D=W^{*}(p_{k}:k\geq 0) \),
\[
L_{k}^{*}dL_{k'}=v^{k}\delta _{kk'}\theta ((v^{*})^{k}dv^{k})(v^{*})^{k},\]
with the convention that \( v^{-1}=v^{*} \). 

By \cite{speicher:thesis} and \cite[Example 2.6]{shlyakht:semicirc}, \( L_{k} \)
are free with amalgamation over \( D \). 

Let \( M_{1}=W^{*}(D,L_{k}+L_{k}^{*},k\leq 0) \).  Since when \( k\leq 0 \),
\[
L_{k}^{*}dL_{k}=C_{k}\theta (d),\]
for some nonzero constants \( C_{k} \), we get by \cite[Theorem 2.3]{shlyakht:amalg}
that \( L_{k} \) is \( * \)-free from \( D \), and
\[
M_{1}=(D,\theta )*(*_{k\leq 0}W^{*}(L_{k}+L_{k}^{*})).\]
Since \( L_{k}+L_{k}^{*} \) is semicircular, we get
\[
M_{1}=(D,\theta )*L(\mathbb {F}_{\infty })\cong L(\mathbb {F}_{\infty }),\]
by  the results of Ken Dykema \cite{dykema:interpolated}. 

Now, \( M=M_{1}*_{D}(W^{*}(L_{k}+L_{k}^{*},k>0)). \)  Since \( M_{1} \) is
a factor and \( L_{k} \) are \( \omega  \)-creation  operators over \( M \),
where \( \omega _{ij}=\delta _{ij}C_{i}p_{i}\theta ((v^{*})^{i}dv^{i}), \)
we get by \cite[Proposition 5.4]{shlyakht:amalg} that   \( M\cong M_{1}*L(\mathbb {F}_{t}), \)
where \( t=\sum _{i}\theta (v^{i}(v^{*})^{i})^{2} \).  Hence \( M\cong L(\mathbb {F}_{\infty } \)). 
\end{proof}

\section{Quantum \protect\( SU_{q}\protect \) actions.\label{sec:quantumactions}}

The main goal of this section is to analyze inclusions of type III factors coming
from a minimal free product type action of \( SU_{q}(2) \) on a von Neumann
algebra. 

Let \( (A,\delta ) \) be the Hopf-von~Neumann algebra of Woronowicz's quantum
group \( SU_{q}(2) \) \cite{woronowicz:SUq}, \cite{woronowicz:compactMatrixPseudogroups},
and denote by \( h \) the canonical Haar state on \( A \). Let \( V \) be
the multiplicative unitary, i.e., the unitary on \( L^{2}(A)\otimes L^{2}(A) \)
defined by
\[
V(\Lambda _{h}(a)\otimes \zeta )=\delta (a)(\Lambda _{h}(1)\otimes \zeta ),\]
and let \( W \) be the fundamental unitary, i.e., unitary on \( L^{2}(A)\otimes L^{2}(A) \)
defined by
\[
W(\eta \otimes \Lambda _{h}(a))=\delta (a)(\eta \otimes \Lambda _{h}(1)),\]
where \( \Lambda _{h}:A\rightarrow L^{2}(A) \) is the canonical injection associated
with the Haar state \( h \). The dual Hopf-von Neumann algebra \( \hat{A} \)
is the von Neumann algebra acting on \( L^{2}(A) \) generated by elements of
the form:
\[
(\text{Id}\otimes \omega )(V),\quad \omega \in B(L^{2}(A))_{*}.\]
Let \( m \) be an integer \( \geqq 2 \) or \( +\infty  \) and set
\[
M=(L(\mathbb {F}_{m}),\tau )*(A,h)\]
and let \( \Gamma :M\to M\otimes A \) be the free product of the trivial action
of \( SU_{q}(2) \) on \( L(\mathbb {F}_{m}) \), and the left regular representation
action \( \delta  \) of \( SU_{q}(2) \) on \( A \) \cite[\S3]{ueda:fixed-point}. 

Recall that the crossed-product \( M\rtimes _{\Gamma }SU_{q}(2) \) is the von
Neumann subalgebra of \( M\otimes B(L^{2}(A)) \) generated by \( \Gamma (M) \)
and \( {\mathbb C}1\otimes \hat{A} \).

For each irreducible finite-dimensional unitary representation \( \pi :V_{\pi }\to V_{\pi }\otimes A \)
of \( SU_{q}(2) \), consider the Wassermann-type inclusion \cite{wassermann:ergodicActionsAndCoactionsSubfactors}
of von Neumann algebras
\[
\mathcal{N}=\mathbb {C}1\otimes M^{\Gamma }\subset (B(V_{\pi })\otimes M)^{\textrm{Ad}{\pi }\otimes \Gamma }=\mathcal{M}.\]
It turns out that the negative (normalized) \( q \)-trace \( \tau _{q}^{(\pi ,-)} \)
gives rise to a finite-index conditional expectation \( \mathcal{E}:\mathcal{M}\to \mathcal{N} \)
(see \cite[eq. (3.6)]{ueda:fixed-point}). We collect below some facts about
this inclusion, most of which are from \cite{ueda:fixed-point}.

\begin{thm}
\label{thrm:propertiescrossedproduct} Let \( \pi  \) be an irreducible representation
of \( SU_{q}(2) \) and let \( \mathcal{N}\subset \mathcal{M} \) be as above,
taken with the conditional expectation \( \mathcal{E} \). Then\\
(a) \( \mathcal{N} \) and \( \mathcal{M} \) are type III\( _{q^{2}} \) factors
and \( \mathcal{N}\cong \mathcal{M} \);\\
(b) The inclusion \emph{\( \mathcal{N}\subset \mathcal{M} \)} is irreducible;\\
(c) The index of \( \mathcal{E} \) is the square of the \( q \)-dimension
of \( \pi  \), \( (\dim _{q}\pi )^{2} \);\\
(d) Let \( \phi  \) denote the restriction of the free product state \( \tau *h \)
to \( \mathcal{N}=M^{\Gamma }\subset M \). Then the centralizers \( \hat{\mathcal{N}}=\mathcal{N}^{\phi } \)
and \( \hat{\mathcal{M}}=\mathcal{M}^{\phi \circ \mathcal{E}} \) are both factors
of type II\( _{1} \);\\
(e) The inclusion \( \mathcal{N}\subset \mathcal{M} \) is ``essentially type
II''. More precisely, the restriction of \( \mathcal{E} \) to the inclusion
\( N\subset M \) is trace-preserving, and
\[
\hat{\mathcal{N}}=\mathcal{N}^{\phi }\subset \mathcal{M}^{\phi \circ \mathcal{E}}=\hat{\mathcal{M}}\]
is an inclusion of type II\( _{1} \) factors with the same index and the same
system of higher relative commutants as \( \mathcal{N}\subset \mathcal{M} \).
In particular, \( \hat{\mathcal{N}}\subset \hat{\mathcal{M}} \) is irreducible.\\
(f) The algebras \( \mathcal{N}=M^{\Gamma } \) and \( M\rtimes _{\Gamma }SU_{q}(2) \)
are isomorphic.
\end{thm}
\begin{proof}
The fact that \( \mathcal{N} \) is of type III\( _{q^{2}} \) was proved in \cite{ueda:fixed-point};
this fact also follows from the explicit description of \( \mathcal{N} \) given
in Appendix I.

For the convenience of the reader we give an expanded proof of the isomorphism
\( \mathcal{N}\cong \mathcal{M} \) (see \cite[Remark 10]{ueda:fixed-point}).
Let \( n \) be the dimension of \( V_{\pi } \) as a vector space. Choose a
standard orthonormal system of vectors \( \xi _{i} \), \( 1\leq i\leq n \)
for the space \( V_{\pi } \) so that the (co)representation \( \pi  \) is
given by
\[
\pi (\xi _{j})=\sum _{i=1}^{n}\xi _{i}\otimes u_{ij}\in V_{\pi }\otimes A,\]
for some \( u_{ij}\in A \) satisfying \( \sum _{k}u_{ki}^{*}u_{kj}=\delta _{ij} \). 

Since \( \mathcal{N}=M^{\Gamma }\subset M \) is of type III, we can find isometries
\( S_{i}\in M^{\Gamma } \), \( i=1,\dots ,n \) so that \( \sum S_{i}S_{i}^{*}=1 \)
and \( S_{i}^{*}S_{j}=\delta _{ij} \). Set
\[
v_{j}=\sum _{i=1}^{n}S_{i}\cdot u_{ij}\in M^{\Gamma }\cdot A\subset M,\quad H=\textrm{span}_{\mathbb {C}}\{v_{j}:1\leq j\leq 
n\}.\]
Then \( v_{i}^{*}v_{j}=\delta _{ij} \) and \( \sum v_{i}v_{i}^{*}=1 \). Moreover,
since \( S_{i}\in M^{\Gamma } \) are fixed by \( \Gamma  \), we get (since
\( \Gamma  \) is an action and \( \Gamma |_{A}=\delta  \))
\begin{eqnarray*}
\Gamma (v_{j}) & = & \sum _{i}(S_{j}\otimes 1)\cdot \Gamma (u_{ij})\\
 & = & \sum _{i}(S_{j}\otimes 1)\cdot \delta (u_{ij})\\
 & = & \sum _{i,k}(S_{j}\otimes 1)\cdot (u_{ki}\otimes u_{ij})\\
 & = & \sum _{i}V_{i}\otimes u_{ij}.
\end{eqnarray*}
It follows that the isomorphism \( H\ni v_{i}\mapsto \xi _{i}\in V_{\pi } \)
is equivariant with respect to the restriction of \( \Gamma  \) to \( H \)
and the (co)representation \( \pi  \) on \( V \). It follows that the isomorphism
\[
\Phi :B(V_{\pi })\otimes M\to M,\quad \Phi ((q_{ij})_{ij})=\sum _{ij}v_{i}q_{ij}v_{j}^{*}\]
is equivariant with respect to the action \( \textrm{Ad}\pi \otimes \Gamma  \)
on \( B(V_{\pi })\otimes M \) and the action \( \Gamma  \) on \( M \). Thus
\( \Phi  \) restricts to an isomorphism of \( \mathcal{N}=M^{\Gamma } \) and
\( \mathcal{M}=(B(V_{\pi })\otimes M)^{(\textrm{Ad}\pi \otimes \Gamma )} \).
Thus (a) holds.

Statements (b)--(e) were proved in \cite{ueda:fixed-point} (Theorem 8 and the
discussion on pp. 43--49).

For part (f), it is known that
\[
M\rtimes _{\Gamma }SU_{q}(2)=(M\otimes B(L^{2}(A)))^{\widetilde{\Gamma }},\]
where \( \widetilde{\Gamma }:=\text {Ad}(1\otimes \Sigma (W^{*}))\circ (\text {Id}\otimes \Sigma )\circ (\Gamma \otimes \text 
{Id}) \),
and where \( \Sigma  \) is the flip map (see \cite[Remark 20]{boca:ergodicActrionsQuantumGroups}
for a simple proof, which also works in the \( W^{*} \)-algebraic setting without
any change. For the reader's convenience, we should mention that the definition
of ``\( \lambda (\omega ) \)'' in that paper involves a typographical error.)
Since \( (M,\Gamma ) \) contains \( (A,\delta ) \), we see that
\[
(M\otimes B(L^{2}(A)))^{\widetilde{\Gamma }}\cong M^{\Gamma }\otimes B(L^{2}(A))\]
(see \cite[Lemma 4.2]{ueda:SUqActionFullFactor}) so that the crossed-product
\( M\rtimes _{\Gamma }SU_{q}(2) \) and the fixed-point algebra \( M^{\Gamma } \)
are isomorphic to each other, both being of type III. 
\end{proof}
We would like to point out that the \( n \)-th step \( \mathcal{M}_{n} \)
in the basic construction associated to \( \mathcal{N}\subset \mathcal{M} \)
is rather easy to describe:
\[
\mathcal{M}_{n}\cong (B(V_{\rho _{n}})\otimes M)^{\textrm{Ad}\rho _{n}\otimes \Gamma }\]
where
\begin{eqnarray*}
\rho _{2s+1} & = & \hat{\pi }\otimes \rho _{2s}\\
\rho _{2s} & = & \pi \otimes \rho _{2s-1}\\
\rho _{0} & = & \pi ,
\end{eqnarray*}
and \( \hat{\pi } \) is the canonical dual of \( \pi  \). The conditional
expectation \( \mathcal{E}_{n}:\mathcal{M}_{n}\to \mathcal{M}_{n-1} \) is given
by
\begin{eqnarray*}
\mathcal{E}_{2s+1} & = & \tau _{q}^{(\hat{\pi },-)}\otimes \textrm{Id}\otimes \cdots \\
\mathcal{E}_{2s} & = & \tau _{q}^{(\pi ,-)}\otimes \textrm{Id}\otimes \cdots .
\end{eqnarray*}
This allows one to compute the \( \lambda  \)-lattice of higher relative commutants
of \( \mathcal{N}\subset \mathcal{M} \) (and thus of \( \hat{\mathcal{N}}\subset \hat{\mathcal{M}} \))
in terms of the representation theory of \( SU_{q}(2) \). In particular, one
can show that if \( \pi  \) is taken to be the fundamental representation,
then the principal graph of \( \hat{\mathcal{N}}\subset \hat{\mathcal{M}} \)
is \( A_{\infty } \), and the index is \( (q+q^{-1})^{2} \).

For \( 0<\lambda <1 \), let \( T_{\lambda } \) denote the (unique) free Araki-Woods
factor of type III\( _{\lambda } \) (see \cite{shlyakht:quasifree:big} and
Section \ref{sec:amalgfreeprod}).

\begin{prop}
\label{prop:crossedproductisomTlambda}The factor \( ((L(\mathbb {F}_m),\tau )*(A,h))\rtimes _{\Gamma }SU_{q}(2) \)
is isomorphic to the free Araki-Woods factor \( T_{q^{2}} \). Hence \( \mathcal{M}\cong \mathcal{N}\cong T_{q^{2}} \).
\end{prop}
\begin{proof}
It is known that the fundamental unitary \( W \) lies in \( A\otimes \hat{A}' \)
(acting on \( L^{2}(A)\otimes L^{2}(A) \)), and hence the mapping \( \Phi :x\in B(L^{2}(A))\mapsto W(1\otimes x)W^{*} \)
gives rise to an isomorphism from \( B(L^{2}(A)) \) onto \( A\rtimes _{\delta }SU_{q}(2)=A\rtimes _{\delta }\hat{A} \)
(see e.g. \cite{baaj-skandalis:multiplicativeUnitaries}). 

Moreover, the crossed product structure and the Haar state \( h \) on \( A \)
give rise to a normal faithful conditional expectation
\[
E:B(L^{2}(A))=A\rtimes _{\delta }\hat{A}\to \hat{A}\]
given by
\[
E(x)=(h\otimes \text {\textrm{Id}})(W(1\otimes x)W^{*}),\quad x\in B(L^{2}(A)).\]
We claim that
\[
M\rtimes _{\Gamma }SU_{q}(2)\cong (L(\mathbb {F}_{m})\otimes \hat{A},\tau \otimes \id )*_{\hat{A}}(B(L^{2}(A)),E).\]

This isomorphism is obtained by extending the map
\[
\Psi := \Phi^{-1} :A\rtimes _{\delta }SU_{q}(2)\to B(L^{2}(A))\]
to
\[
[(L(\mathbb {F}_{m}),\tau )*(A,h)]\rtimes _{\Gamma }SU_{q}(2)\supset A\rtimes _{\delta }SU_{q}(2)\]
by mapping elements in \( L(\mathbb {F}_{m}) \) to the canonical copy of \( L(\mathbb {F}_{m}) \)
inside \( (L(\mathbb {F}_{m})\otimes \hat{A},\tau \otimes \id )*_{\hat{A}}(B(L^{2}(A)),E) \).
Since \( \Phi (\hat{a})=1\otimes \hat{a} \), \( \hat{a}\in \hat{A} \) and
\( (h\otimes \text {Id})\circ \Phi =E \), \( \Psi = \Phi^{-1} \) extends to the desired
isomorphism (a freeness condition needs to be checked, see \cite[Proposition 1]{ueda:fixed-point}
and \cite{DVV:book}.) 

(If instead of a quantum group action, we would have an action of an ordinary
group \( G \), the isomorphism above would be the content of the general identity
\[
(R,\phi )*(S,\psi )\rtimes _{\textrm{id}*\alpha \textrm{ }}G=(R\otimes L(G),\phi \otimes \textrm{id})*_{L(G)}(S\rtimes _{\alpha 
}G,E_{G}),\]
where \( E_{G} \) is the canonical \( L(G) \)-valued conditional expectation
on \( S\rtimes _{\alpha }G \) associated to the state \( \psi :S\to \mathbb {C} \).)

Applying Theorem \ref{thrm:freeproductBofH}, and noting that \( M\rtimes _{\Gamma }SU_{q}(2) \)
is of type III\( _{q^{2}} \) (Theorem \ref{thrm:propertiescrossedproduct} (f)
and (a)), we find that
\[
\mathcal{M}\cong \mathcal{N}\cong (L(\mathbb {F}_{m})\otimes \hat{A},\tau \otimes \id )*_{\hat{A}}(B(H),E)\]
is isomorphic to the unique type III\( _{q^{2}} \) free Araki-Woods factor \( T_{q^{2}} \).
\end{proof}
\begin{rem}
Since the free Araki-Woods factor \( T_{q^{2}} \) is prime (i.e., cannot be
written as a tensor product of two diffuse von Neumann algebras, see \cite{shlyakht:prime}),
it follows
that the inclusion \( \mathcal{N}\subset \mathcal{M} \) cannot be decomposed
as a tensor product of a type II inclusion with a fixed type III factor, although
the \( \lambda  \)-lattices of the type III and type II inclusions coincide.
(This remark can be also obtained from the results of Ge \cite{ge:entropy2}).
\end{rem}
\begin{lem}
\label{lemma:centralizerLFinfty}With the assumptions of Theorem \ref{thrm:propertiescrossedproduct}
(d), \( \hat{\mathcal{M}} \) and \( \hat{\mathcal{N}} \) are both isomorphic
to \( L(\mathbb {F}_{\infty }) \). 
\end{lem}
\begin{proof}
Since \( \hat{N}=\mathcal{N}^{\phi } \) is a factor (Theorem \ref{thrm:propertiescrossedproduct}
(d)) and \( \mathcal{N}=M^{\Gamma } \) is isomorphic to \( T_{q^{2}} \) ,
\( \hat{\mathcal{N}}\cong L(\mathbb {F}_{\infty }) \) (see Theorem \ref{thrm:coreofTlambda}).
Similarly, since \( \mathcal{M}\cong T_{q^{2}} \) and \emph{\( \hat{\mathcal{M}}=\mathcal{M}^{\phi \circ \mathcal{E}} \)}
is a factor, \( \hat{\mathcal{M}}\cong L(\mathbb {F}_{\infty }) \).
\end{proof}
For a general finite dimensional representation \( \pi  \) the factor map structure
associated with the type III inclusion \( \mathcal{M}\supseteq \mathcal{N} \)
(see \cite{kosaki:typeIIIsubfactors}) can be described as follows. Decompose
\[
\pi =n_{1}\cdot \pi _{\ell _{1}}\oplus \cdots \oplus n_{k}\cdot \pi _{\ell _{k}}\]
into multiples of irreducible spin \( \ell  \) representations \( \pi _{\ell } \)
(see \cite{masudaANDCO:representationsSUQ2andQJacobiPolynomials}) with \( \ell _{i}\neq \ell _{j} \)
(\( i\neq j \)). 

Set \( X:=\{1,2,\dots ,k\}\times [0,-\log q^{2}) \) and \( F_{t}(j,s):=(j,s+t) \),
and let \( (X_{\mathcal{M}},F_{t}^{\mathcal{M}}) \) and \( (X_{\mathcal{N}},F_{t}^{\mathcal{N}}) \)
be the flows of weights associated with \( \mathcal{M} \) and \( \mathcal{N} \),
respectively, which are both identified with \( ([0,-\log q^{2}),\textrm{translation by }t) \).
Let us define two factor maps
\[
\pi^{\mathcal{M}}:X\rightarrow X_{\mathcal{M}}=[0,-\log q^{2}),\quad \pi ^{\mathcal{N}}:X\rightarrow X_{\mathcal{N}}=[0,-\log 
q^{2})\]
by
\[
\pi^{\mathcal{M}}(j,s)=s-\ell _{j}\log q^{2}\mod -\log q^{2},\quad \pi ^{\mathcal{N}}(j,s)=s.\]
These \( \pi ^{\mathcal{M}} \), \( \pi ^{\mathcal{N}} \) describe the factor
map structure. As a consequence, we get: 

\begin{rem}
The inclusion \( \mathcal{M}\supseteq \mathcal{N} \) is of essentially type
II if and only if all the set
\[
\{\ell :\pi _{\ell }\textrm{ occurs in the decomposition of }\pi \}\]
either consists entirely of half-integers or consists entirely of integers.
Under this assumption, all statements in Theorem \ref{thrm:propertiescrossedproduct}
still hold, with the exception of statement (b) (the resulting inclusion is
of course no longer irreducible).
\end{rem}
This means that our method of constructing irreducible type II\( _{1} \) subfactors
can be applied to the special standard \( \lambda  \)-lattices arising from
representations of the quantum group \( SU_{q}(2) \) (see \cite{banica:compactQuantumGroupsSubfactors}).

We are now ready to prove our main result:

\begin{thm}
For every \( \lambda \in \mathcal{I} \), there exists an irreducible subfactor
\( N \) of \( M=L(\mathbb {F}_{\infty }) \) of index \( \lambda  \), so that
\( N\cong L(\mathbb {F}_{\infty }) \).
\end{thm}
\begin{proof}
The case \( \lambda \in \mathcal{I}\cap [0,4] \) was obtained by Radulescu
in \cite{radulescu:subfact}.

Let \( \lambda >4 \), and choose \( q \) so that \( \lambda =(q+q^{-1})^{2} \).
Let \( \pi  \) be the fundamental representation of \( SU_{q}(2) \). By Theorem
\ref{thrm:propertiescrossedproduct} (e) and Lemma \ref{lemma:centralizerLFinfty},
for each \( 0<q<1 \), we obtain an irreducible inclusion of type II\( _{1} \)
factors
\[
L(\mathbb {F}_{\infty })\cong \hat{\mathcal{N}}=\mathcal{N}^{\phi }\subset \mathcal{M}^{\phi \circ 
\mathcal{E}}=\hat{\mathcal{M}}\cong L(\mathbb {F}_{\infty })\]
of index \( (q+q^{-1})^{2} \). 
\end{proof}
\begin{rem}
In the same way, we see that any irreducible representation \( \pi  \) of \( SU_{q}(2) \)
gives rise to an irreducible subfactor \( \hat{\mathcal{N}}\subset \hat{\mathcal{M}} \),
\( \hat{\mathcal{N}}\cong \hat{\mathcal{M}}\cong L(\mathbb {F}_{\infty }) \),
of index \( (\dim _{q}\pi )^{2} \).
\end{rem}

\section*{Appendix I.}

We give here a concrete system of generators for the von Neumann algebra \( M^{\Gamma } \)
appearing in the proof of the main theorem. This gives also a concrete example
illustrating the proof of Theorem \ref{thrm:genresultA-valued} works.

Let \( H \) be a Hilbert space, given as a direct sum
\[
H=\bigoplus _{k}H_{k}\otimes K_{k},\]
with \( H_{k} \) and \( K_{k} \) finite-dimensional or infinite-dimensional.
Let \( \hat{A} \) be the subalgebra of \( B(H) \) given by
\[
\bigoplus _{k}B(H_{k})\otimes \textrm{Id}_{K_{k}}.\]
For each \( k \), let \( \theta _{k} \) be a normal faithful state on \( B(K_{k}) \)
given by
\[
\theta _{k}(T)=\textrm{Tr}(D_{k}T),\]
where \( D_{k}\in B(K_{k}) \) is a fixed positive matrix. Then
\[
E=\bigoplus _{k}\id \otimes \theta _{k}(\cdot )1_{B(K_{k})}\]
is a normal faithful conditional expectation from \( B(H) \) onto \( \hat{A} \).

Let \( d_{k}=\dim K_{k} \). Choose a basis for \( K_{k} \), so that \( D_{k} \)
is diagonal, with eigenvalues \( \lambda _{1}(k),\ldots ,\lambda _{d(k)}(k) \).
Write \( e_{ij}(k) \), \( 1\leq i,j\leq d(k) \) for the matrix units generating
\( B(K_{k}) \), associated to this basis.

Let \( l^{k}_{ij}(p) \), \( p,k=1,2,\ldots  \), \( 1\leq i,j\leq d(k) \)
be isometries satisfying \( l_{ij}^{k}(p)\cdot l_{i'j'}^{k'}(p')=\delta _{pp'}\delta _{ii'}\delta _{jj'}\delta _{kk'}1 \),
and generating the Toeplitz extension \( \mathcal{T} \) of the Cuntz algebra.
Let \( \psi  \) be the canonical vacuum state on \( \mathcal{T} \) (see \cite{DVV:book}).
Consider the tensor product conditional expectation
\( \psi \otimes E  \) on \( P=B(L^2(\mathcal{T},\psi)) \otimes B(H) \).

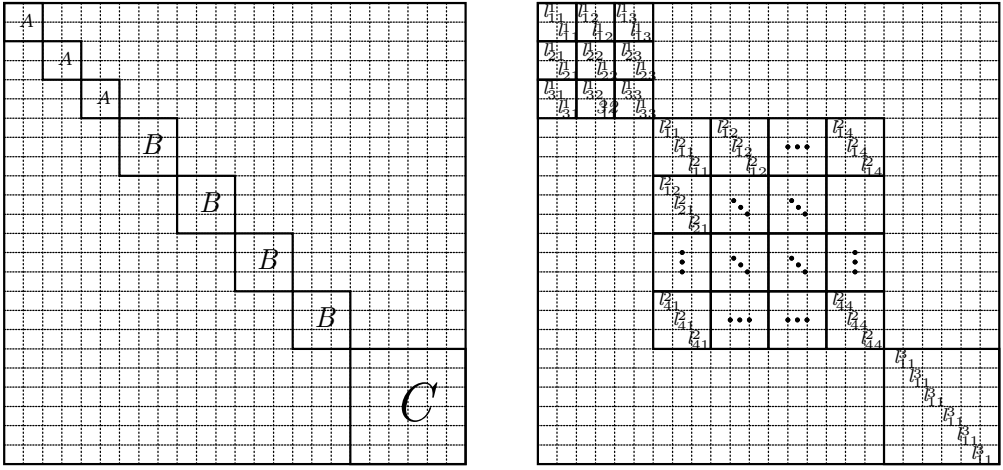
\begin{figure} 

\setlength{\unitlength}{0.00033333in}
\begingroup\makeatletter\ifx\SetFigFont\undefined%
\gdef\SetFigFont#1#2#3#4#5{%
  \reset@font\fontsize{#1}{#2pt}%
  \fontfamily{#3}\fontseries{#4}\fontshape{#5}%
  \selectfont}%
\fi\endgroup%
{
\begin{picture}(15569,7259)(0,-10)
\put(8422,6997){\makebox(0,0)[lb]{\smash{{{\SetFigFont{7}{8.4}{\rmdefault}{\mddefault}{\itdefault}l}}}}}
\put(8497,7072){\makebox(0,0)[lb]{\smash{{{\SetFigFont{5}{6.0}{\rmdefault}{\mddefault}{\updefault}1}}}}}
\put(8497,6922){\makebox(0,0)[lb]{\smash{{{\SetFigFont{5}{6.0}{\rmdefault}{\mddefault}{\updefault}11}}}}}
\put(8647,6697){\makebox(0,0)[lb]{\smash{{{\SetFigFont{7}{8.4}{\rmdefault}{\mddefault}{\itdefault}l}}}}}
\put(8722,6772){\makebox(0,0)[lb]{\smash{{{\SetFigFont{5}{6.0}{\rmdefault}{\mddefault}{\updefault}1}}}}}
\put(8722,6622){\makebox(0,0)[lb]{\smash{{{\SetFigFont{5}{6.0}{\rmdefault}{\mddefault}{\updefault}11}}}}}
\put(10222,5197){\makebox(0,0)[lb]{\smash{{{\SetFigFont{7}{8.4}{\rmdefault}{\mddefault}{\itdefault}l}}}}}
\put(10297,5272){\makebox(0,0)[lb]{\smash{{{\SetFigFont{5}{6.0}{\rmdefault}{\mddefault}{\updefault}2}}}}}
\put(10297,5122){\makebox(0,0)[lb]{\smash{{{\SetFigFont{5}{6.0}{\rmdefault}{\mddefault}{\updefault}11}}}}}
\put(10447,4897){\makebox(0,0)[lb]{\smash{{{\SetFigFont{7}{8.4}{\rmdefault}{\mddefault}{\itdefault}l}}}}}
\put(10522,4972){\makebox(0,0)[lb]{\smash{{{\SetFigFont{5}{6.0}{\rmdefault}{\mddefault}{\updefault}2}}}}}
\put(10522,4822){\makebox(0,0)[lb]{\smash{{{\SetFigFont{5}{6.0}{\rmdefault}{\mddefault}{\updefault}11}}}}}
\put(10672,4597){\makebox(0,0)[lb]{\smash{{{\SetFigFont{7}{8.4}{\rmdefault}{\mddefault}{\itdefault}l}}}}}
\put(10747,4672){\makebox(0,0)[lb]{\smash{{{\SetFigFont{5}{6.0}{\rmdefault}{\mddefault}{\updefault}2}}}}}
\put(10747,4522){\makebox(0,0)[lb]{\smash{{{\SetFigFont{5}{6.0}{\rmdefault}{\mddefault}{\updefault}11}}}}}
\put(13897,1597){\makebox(0,0)[lb]{\smash{{{\SetFigFont{7}{8.4}{\rmdefault}{\mddefault}{\itdefault}l}}}}}
\put(13972,1672){\makebox(0,0)[lb]{\smash{{{\SetFigFont{5}{6.0}{\rmdefault}{\mddefault}{\updefault}3}}}}}
\put(13972,1522){\makebox(0,0)[lb]{\smash{{{\SetFigFont{5}{6.0}{\rmdefault}{\mddefault}{\updefault}11}}}}}
\put(14122,1297){\makebox(0,0)[lb]{\smash{{{\SetFigFont{7}{8.4}{\rmdefault}{\mddefault}{\itdefault}l}}}}}
\put(14197,1372){\makebox(0,0)[lb]{\smash{{{\SetFigFont{5}{6.0}{\rmdefault}{\mddefault}{\updefault}3}}}}}
\put(14197,1222){\makebox(0,0)[lb]{\smash{{{\SetFigFont{5}{6.0}{\rmdefault}{\mddefault}{\updefault}11}}}}}
\put(14347,997){\makebox(0,0)[lb]{\smash{{{\SetFigFont{7}{8.4}{\rmdefault}{\mddefault}{\itdefault}l}}}}}
\put(14422,1072){\makebox(0,0)[lb]{\smash{{{\SetFigFont{5}{6.0}{\rmdefault}{\mddefault}{\updefault}3}}}}}
\put(14422,922){\makebox(0,0)[lb]{\smash{{{\SetFigFont{5}{6.0}{\rmdefault}{\mddefault}{\updefault}11}}}}}
\put(14647,697){\makebox(0,0)[lb]{\smash{{{\SetFigFont{7}{8.4}{\rmdefault}{\mddefault}{\itdefault}l}}}}}
\put(14722,772){\makebox(0,0)[lb]{\smash{{{\SetFigFont{5}{6.0}{\rmdefault}{\mddefault}{\updefault}3}}}}}
\put(14722,622){\makebox(0,0)[lb]{\smash{{{\SetFigFont{5}{6.0}{\rmdefault}{\mddefault}{\updefault}11}}}}}
\put(14872,397){\makebox(0,0)[lb]{\smash{{{\SetFigFont{7}{8.4}{\rmdefault}{\mddefault}{\itdefault}l}}}}}
\put(14947,472){\makebox(0,0)[lb]{\smash{{{\SetFigFont{5}{6.0}{\rmdefault}{\mddefault}{\updefault}3}}}}}
\put(14947,322){\makebox(0,0)[lb]{\smash{{{\SetFigFont{5}{6.0}{\rmdefault}{\mddefault}{\updefault}11}}}}}
\put(15097,97){\makebox(0,0)[lb]{\smash{{{\SetFigFont{7}{8.4}{\rmdefault}{\mddefault}{\itdefault}l}}}}}
\put(15172,172){\makebox(0,0)[lb]{\smash{{{\SetFigFont{5}{6.0}{\rmdefault}{\mddefault}{\updefault}3}}}}}
\put(15172,22){\makebox(0,0)[lb]{\smash{{{\SetFigFont{5}{6.0}{\rmdefault}{\mddefault}{\updefault}11}}}}}
\put(12247,4972){\blacken\ellipse{68}{68}}
\put(12247,4972){\ellipse{68}{68}}
\put(12397,4972){\blacken\ellipse{68}{68}}
\put(12397,4972){\ellipse{68}{68}}
\put(12547,4972){\blacken\ellipse{68}{68}}
\put(12547,4972){\ellipse{68}{68}}
\put(11347,2272){\blacken\ellipse{68}{68}}
\put(11347,2272){\ellipse{68}{68}}
\put(11497,2272){\blacken\ellipse{68}{68}}
\put(11497,2272){\ellipse{68}{68}}
\put(11647,2272){\blacken\ellipse{68}{68}}
\put(11647,2272){\ellipse{68}{68}}
\put(12247,2272){\blacken\ellipse{68}{68}}
\put(12247,2272){\ellipse{68}{68}}
\put(12397,2272){\blacken\ellipse{68}{68}}
\put(12397,2272){\ellipse{68}{68}}
\put(12547,2272){\blacken\ellipse{68}{68}}
\put(12547,2272){\ellipse{68}{68}}
\put(10597,3022){\blacken\ellipse{68}{68}}
\put(10597,3022){\ellipse{68}{68}}
\put(10597,3172){\blacken\ellipse{68}{68}}
\put(10597,3172){\ellipse{68}{68}}
\put(10597,3322){\blacken\ellipse{68}{68}}
\put(10597,3322){\ellipse{68}{68}}
\put(13297,3022){\blacken\ellipse{68}{68}}
\put(13297,3022){\ellipse{68}{68}}
\put(13297,3172){\blacken\ellipse{68}{68}}
\put(13297,3172){\ellipse{68}{68}}
\put(13297,3322){\blacken\ellipse{68}{68}}
\put(13297,3322){\ellipse{68}{68}}
\put(11616,3922){\blacken\ellipse{68}{68}}
\put(11616,3922){\ellipse{68}{68}}
\put(11510,4028){\blacken\ellipse{68}{68}}
\put(11510,4028){\ellipse{68}{68}}
\put(11404,4134){\blacken\ellipse{68}{68}}
\put(11404,4134){\ellipse{68}{68}}
\put(12516,3022){\blacken\ellipse{68}{68}}
\put(12516,3022){\ellipse{68}{68}}
\put(12410,3128){\blacken\ellipse{68}{68}}
\put(12410,3128){\ellipse{68}{68}}
\put(12304,3234){\blacken\ellipse{68}{68}}
\put(12304,3234){\ellipse{68}{68}}
\put(12516,3922){\blacken\ellipse{68}{68}}
\put(12516,3922){\ellipse{68}{68}}
\put(12410,4028){\blacken\ellipse{68}{68}}
\put(12410,4028){\ellipse{68}{68}}
\put(12304,4134){\blacken\ellipse{68}{68}}
\put(12304,4134){\ellipse{68}{68}}
\put(11616,3022){\blacken\ellipse{68}{68}}
\put(11616,3022){\ellipse{68}{68}}
\put(11510,3128){\blacken\ellipse{68}{68}}
\put(11510,3128){\ellipse{68}{68}}
\put(11404,3234){\blacken\ellipse{68}{68}}
\put(11404,3234){\ellipse{68}{68}}
\put(11122,5197){\makebox(0,0)[lb]{\smash{{{\SetFigFont{7}{8.4}{\rmdefault}{\mddefault}{\itdefault}l}}}}}
\put(11197,5272){\makebox(0,0)[lb]{\smash{{{\SetFigFont{5}{6.0}{\rmdefault}{\mddefault}{\updefault}2}}}}}
\put(11197,5122){\makebox(0,0)[lb]{\smash{{{\SetFigFont{5}{6.0}{\rmdefault}{\mddefault}{\updefault}12}}}}}
\put(11347,4897){\makebox(0,0)[lb]{\smash{{{\SetFigFont{7}{8.4}{\rmdefault}{\mddefault}{\itdefault}l}}}}}
\put(11422,4972){\makebox(0,0)[lb]{\smash{{{\SetFigFont{5}{6.0}{\rmdefault}{\mddefault}{\updefault}2}}}}}
\put(11422,4822){\makebox(0,0)[lb]{\smash{{{\SetFigFont{5}{6.0}{\rmdefault}{\mddefault}{\updefault}12}}}}}
\put(11572,4597){\makebox(0,0)[lb]{\smash{{{\SetFigFont{7}{8.4}{\rmdefault}{\mddefault}{\itdefault}l}}}}}
\put(11647,4672){\makebox(0,0)[lb]{\smash{{{\SetFigFont{5}{6.0}{\rmdefault}{\mddefault}{\updefault}2}}}}}
\put(11647,4522){\makebox(0,0)[lb]{\smash{{{\SetFigFont{5}{6.0}{\rmdefault}{\mddefault}{\updefault}12}}}}}
\put(10222,4297){\makebox(0,0)[lb]{\smash{{{\SetFigFont{7}{8.4}{\rmdefault}{\mddefault}{\itdefault}l}}}}}
\put(10297,4372){\makebox(0,0)[lb]{\smash{{{\SetFigFont{5}{6.0}{\rmdefault}{\mddefault}{\updefault}2}}}}}
\put(10297,4222){\makebox(0,0)[lb]{\smash{{{\SetFigFont{5}{6.0}{\rmdefault}{\mddefault}{\updefault}12}}}}}
\put(10447,3997){\makebox(0,0)[lb]{\smash{{{\SetFigFont{7}{8.4}{\rmdefault}{\mddefault}{\itdefault}l}}}}}
\put(10522,4072){\makebox(0,0)[lb]{\smash{{{\SetFigFont{5}{6.0}{\rmdefault}{\mddefault}{\updefault}2}}}}}
\put(10522,3922){\makebox(0,0)[lb]{\smash{{{\SetFigFont{5}{6.0}{\rmdefault}{\mddefault}{\updefault}21}}}}}
\put(10672,3697){\makebox(0,0)[lb]{\smash{{{\SetFigFont{7}{8.4}{\rmdefault}{\mddefault}{\itdefault}l}}}}}
\put(10747,3772){\makebox(0,0)[lb]{\smash{{{\SetFigFont{5}{6.0}{\rmdefault}{\mddefault}{\updefault}2}}}}}
\put(10747,3622){\makebox(0,0)[lb]{\smash{{{\SetFigFont{5}{6.0}{\rmdefault}{\mddefault}{\updefault}21}}}}}
\put(10222,2497){\makebox(0,0)[lb]{\smash{{{\SetFigFont{7}{8.4}{\rmdefault}{\mddefault}{\itdefault}l}}}}}
\put(10297,2572){\makebox(0,0)[lb]{\smash{{{\SetFigFont{5}{6.0}{\rmdefault}{\mddefault}{\updefault}2}}}}}
\put(10297,2422){\makebox(0,0)[lb]{\smash{{{\SetFigFont{5}{6.0}{\rmdefault}{\mddefault}{\updefault}41}}}}}
\put(10447,2197){\makebox(0,0)[lb]{\smash{{{\SetFigFont{7}{8.4}{\rmdefault}{\mddefault}{\itdefault}l}}}}}
\put(10522,2272){\makebox(0,0)[lb]{\smash{{{\SetFigFont{5}{6.0}{\rmdefault}{\mddefault}{\updefault}2}}}}}
\put(10522,2122){\makebox(0,0)[lb]{\smash{{{\SetFigFont{5}{6.0}{\rmdefault}{\mddefault}{\updefault}41}}}}}
\put(10672,1897){\makebox(0,0)[lb]{\smash{{{\SetFigFont{7}{8.4}{\rmdefault}{\mddefault}{\itdefault}l}}}}}
\put(10747,1972){\makebox(0,0)[lb]{\smash{{{\SetFigFont{5}{6.0}{\rmdefault}{\mddefault}{\updefault}2}}}}}
\put(10747,1822){\makebox(0,0)[lb]{\smash{{{\SetFigFont{5}{6.0}{\rmdefault}{\mddefault}{\updefault}41}}}}}
\put(12922,5197){\makebox(0,0)[lb]{\smash{{{\SetFigFont{7}{8.4}{\rmdefault}{\mddefault}{\itdefault}l}}}}}
\put(12997,5272){\makebox(0,0)[lb]{\smash{{{\SetFigFont{5}{6.0}{\rmdefault}{\mddefault}{\updefault}2}}}}}
\put(12997,5122){\makebox(0,0)[lb]{\smash{{{\SetFigFont{5}{6.0}{\rmdefault}{\mddefault}{\updefault}14}}}}}
\put(13147,4897){\makebox(0,0)[lb]{\smash{{{\SetFigFont{7}{8.4}{\rmdefault}{\mddefault}{\itdefault}l}}}}}
\put(13222,4972){\makebox(0,0)[lb]{\smash{{{\SetFigFont{5}{6.0}{\rmdefault}{\mddefault}{\updefault}2}}}}}
\put(13222,4822){\makebox(0,0)[lb]{\smash{{{\SetFigFont{5}{6.0}{\rmdefault}{\mddefault}{\updefault}14}}}}}
\put(13372,4597){\makebox(0,0)[lb]{\smash{{{\SetFigFont{7}{8.4}{\rmdefault}{\mddefault}{\itdefault}l}}}}}
\put(13447,4672){\makebox(0,0)[lb]{\smash{{{\SetFigFont{5}{6.0}{\rmdefault}{\mddefault}{\updefault}2}}}}}
\put(13447,4522){\makebox(0,0)[lb]{\smash{{{\SetFigFont{5}{6.0}{\rmdefault}{\mddefault}{\updefault}14}}}}}
\put(12922,2497){\makebox(0,0)[lb]{\smash{{{\SetFigFont{7}{8.4}{\rmdefault}{\mddefault}{\itdefault}l}}}}}
\put(12997,2572){\makebox(0,0)[lb]{\smash{{{\SetFigFont{5}{6.0}{\rmdefault}{\mddefault}{\updefault}2}}}}}
\put(12997,2422){\makebox(0,0)[lb]{\smash{{{\SetFigFont{5}{6.0}{\rmdefault}{\mddefault}{\updefault}44}}}}}
\put(13147,2197){\makebox(0,0)[lb]{\smash{{{\SetFigFont{7}{8.4}{\rmdefault}{\mddefault}{\itdefault}l}}}}}
\put(13222,2272){\makebox(0,0)[lb]{\smash{{{\SetFigFont{5}{6.0}{\rmdefault}{\mddefault}{\updefault}2}}}}}
\put(13222,2122){\makebox(0,0)[lb]{\smash{{{\SetFigFont{5}{6.0}{\rmdefault}{\mddefault}{\updefault}44}}}}}
\put(13372,1897){\makebox(0,0)[lb]{\smash{{{\SetFigFont{7}{8.4}{\rmdefault}{\mddefault}{\itdefault}l}}}}}
\put(13447,1972){\makebox(0,0)[lb]{\smash{{{\SetFigFont{5}{6.0}{\rmdefault}{\mddefault}{\updefault}2}}}}}
\put(13447,1822){\makebox(0,0)[lb]{\smash{{{\SetFigFont{5}{6.0}{\rmdefault}{\mddefault}{\updefault}44}}}}}
\dottedline{45}(8347,6922)(15547,6922)
\dottedline{45}(8347,6622)(15547,6622)
\dottedline{45}(8347,6322)(15547,6322)
\dottedline{45}(8347,5722)(15547,5722)
\dottedline{45}(8347,6022)(15547,6022)
\dottedline{45}(8347,3922)(15547,3922)
\dottedline{45}(8347,4222)(15547,4222)
\dottedline{45}(8347,4522)(15547,4522)
\dottedline{45}(8347,4822)(15547,4822)
\dottedline{45}(8347,5422)(15547,5422)
\dottedline{45}(8347,5122)(15547,5122)
\dottedline{45}(8347,3622)(15547,3622)
\dottedline{45}(8347,3022)(15547,3022)
\dottedline{45}(8347,3322)(15547,3322)
\dottedline{45}(8347,2722)(15547,2722)
\dottedline{45}(8347,2122)(15547,2122)
\dottedline{45}(8347,1522)(15547,1522)
\dottedline{45}(8347,1222)(15547,1222)
\dottedline{45}(8347,622)(15547,622)
\dottedline{45}(8347,322)(15547,322)
\dottedline{45}(8347,922)(15547,922)
\dottedline{45}(8347,2422)(15547,2422)
\dottedline{60}(8647,7222)(8647,22)
\dottedline{60}(8947,7222)(8947,22)
\dottedline{60}(9247,7222)(9247,22)
\dottedline{60}(9847,7222)(9847,22)
\dottedline{60}(10147,7222)(10147,22)
\dottedline{60}(10747,7222)(10747,22)
\dottedline{60}(11047,7222)(11047,22)
\dottedline{60}(11347,7222)(11347,22)
\dottedline{60}(11647,7222)(11647,22)
\dottedline{60}(11947,7222)(11947,22)
\dottedline{60}(12547,7222)(12547,22)
\dottedline{60}(12847,7222)(12847,22)
\dottedline{60}(13147,7222)(13147,22)
\dottedline{60}(13447,7222)(13447,22)
\dottedline{60}(14047,7222)(14047,22)
\dottedline{60}(14347,7222)(14347,22)
\dottedline{60}(14647,7222)(14647,22)
\dottedline{60}(14947,7222)(14947,22)
\dottedline{60}(15247,7222)(15247,22)
\dottedline{60}(12247,7222)(12247,22)
\dottedline{60}(13747,7222)(13747,22)
\dottedline{60}(10447,7222)(10447,22)
\dottedline{60}(9547,7222)(9547,22)
\thicklines
\path(8347,7222)(15547,7222)(15547,22)
	(8347,22)(8347,7222)
\thinlines
\dottedline{45}(22,1822)(7222,1822)
\dottedline{45}(8347,1822)(15547,1822)
\dottedline{45}(22,6922)(7222,6922)
\dottedline{45}(22,6322)(7222,6322)
\dottedline{45}(22,5722)(7222,5722)
\dottedline{45}(22,6022)(7222,6022)
\dottedline{45}(22,3922)(7222,3922)
\dottedline{45}(22,4222)(7222,4222)
\dottedline{45}(22,4522)(7222,4522)
\dottedline{45}(22,4822)(7222,4822)
\dottedline{45}(22,5422)(7222,5422)
\dottedline{45}(22,5122)(7222,5122)
\dottedline{45}(22,3622)(7222,3622)
\dottedline{45}(22,3022)(7222,3022)
\dottedline{45}(22,3322)(7222,3322)
\dottedline{45}(22,2722)(7222,2722)
\dottedline{45}(22,2122)(7222,2122)
\dottedline{45}(22,1522)(7222,1522)
\dottedline{45}(22,1222)(7222,1222)
\dottedline{45}(22,622)(7222,622)
\dottedline{45}(22,322)(7222,322)
\dottedline{45}(22,922)(7222,922)
\dottedline{45}(22,2422)(7222,2422)
\dottedline{60}(922,7222)(922,22)
\dottedline{60}(1522,7222)(1522,22)
\dottedline{60}(1822,7222)(1822,22)
\dottedline{60}(2422,7222)(2422,22)
\dottedline{60}(2722,7222)(2722,22)
\dottedline{60}(3022,7222)(3022,22)
\dottedline{60}(3322,7222)(3322,22)
\dottedline{60}(3622,7222)(3622,22)
\dottedline{60}(4222,7222)(4222,22)
\dottedline{60}(4522,7222)(4522,22)
\dottedline{60}(4822,7222)(4822,22)
\dottedline{60}(5122,7222)(5122,22)
\dottedline{60}(5722,7222)(5722,22)
\dottedline{60}(6022,7222)(6022,22)
\dottedline{60}(6322,7222)(6322,22)
\dottedline{60}(6622,7222)(6622,22)
\dottedline{60}(6922,7222)(6922,22)
\dottedline{60}(3922,7222)(3922,22)
\dottedline{60}(5422,7222)(5422,22)
\dottedline{60}(2122,7222)(2122,22)
\dottedline{60}(1222,7222)(1222,22)
\thicklines
\path(622,6622)(1222,6622)(1222,6022)
	(622,6022)(622,6622)
\path(1222,6022)(1822,6022)(1822,5422)
	(1222,5422)(1222,6022)
\path(1822,5422)(2722,5422)(2722,4522)
	(1822,4522)(1822,5422)
\path(2722,4522)(3622,4522)(3622,3622)
	(2722,3622)(2722,4522)
\path(3622,3622)(4522,3622)(4522,2722)
	(3622,2722)(3622,3622)
\path(4522,2722)(5422,2722)(5422,1822)
	(4522,1822)(4522,2722)
\path(5422,1822)(7222,1822)(7222,22)
	(5422,22)(5422,1822)
\path(8947,6622)(9547,6622)(9547,6022)
	(8947,6022)(8947,6622)
\path(9547,6022)(10147,6022)(10147,5422)
	(9547,5422)(9547,6022)
\path(11047,4522)(11947,4522)(11947,3622)
	(11047,3622)(11047,4522)
\path(11947,3622)(12847,3622)(12847,2722)
	(11947,2722)(11947,3622)
\path(12847,2722)(13747,2722)(13747,1822)
	(12847,1822)(12847,2722)
\path(13747,1822)(15547,1822)(15547,22)
	(13747,22)(13747,1822)
\path(22,7222)(7222,7222)(7222,22)
	(22,22)(22,7222)
\thinlines
\dottedline{60}(322,7222)(322,22)
\dottedline{60}(622,7222)(622,22)
\thicklines
\path(22,7222)(622,7222)(622,6622)
	(22,6622)(22,7222)
\thinlines
\dottedline{45}(22,6622)(7222,6622)
\thicklines
\path(9547,6622)(10147,6622)(10147,6022)
	(9547,6022)(9547,6622)
\path(9547,7222)(10147,7222)(10147,6622)
	(9547,6622)(9547,7222)
\path(8947,7222)(9547,7222)(9547,6622)
	(8947,6622)(8947,7222)
\path(8347,6622)(8947,6622)(8947,6022)
	(8347,6022)(8347,6622)
\path(8347,6022)(8947,6022)(8947,5422)
	(8347,5422)(8347,6022)
\path(8947,6022)(9547,6022)(9547,5422)
	(8947,5422)(8947,6022)
\path(8347,7222)(8947,7222)(8947,6622)
	(8347,6622)(8347,7222)
\path(11947,5422)(12847,5422)(12847,4522)
	(11947,4522)(11947,5422)
\path(12847,5422)(13747,5422)(13747,4522)
	(12847,4522)(12847,5422)
\path(12847,4522)(13747,4522)(13747,3622)
	(12847,3622)(12847,4522)
\path(11947,4522)(12847,4522)(12847,3622)
	(11947,3622)(11947,4522)
\path(12847,3622)(13747,3622)(13747,2722)
	(12847,2722)(12847,3622)
\path(11947,2722)(12847,2722)(12847,1822)
	(11947,1822)(11947,2722)
\path(11047,3622)(11947,3622)(11947,2722)
	(11047,2722)(11047,3622)
\path(10147,3622)(11047,3622)(11047,2722)
	(10147,2722)(10147,3622)
\path(11047,2722)(11947,2722)(11947,1822)
	(11047,1822)(11047,2722)
\path(10147,2722)(11047,2722)(11047,1822)
	(10147,1822)(10147,2722)
\path(11047,5422)(11947,5422)(11947,4522)
	(11047,4522)(11047,5422)
\path(10147,5422)(11047,5422)(11047,4522)
	(10147,4522)(10147,5422)
\path(10147,4522)(11047,4522)(11047,3622)
	(10147,3622)(10147,4522)
\put(847,6247){\makebox(0,0)[lb]{\smash{{{\SetFigFont{7}{8.4}{\rmdefault}{\mddefault}{\itdefault}A}}}}}
\put(1447,5647){\makebox(0,0)[lb]{\smash{{{\SetFigFont{7}{8.4}{\rmdefault}{\mddefault}{\itdefault}A}}}}}
\put(2159,4859){\makebox(0,0)[lb]{\smash{{{\SetFigFont{11}{13.2}{\rmdefault}{\mddefault}{\itdefault}B}}}}}
\put(3059,3959){\makebox(0,0)[lb]{\smash{{{\SetFigFont{11}{13.2}{\rmdefault}{\mddefault}{\itdefault}B}}}}}
\put(3959,3059){\makebox(0,0)[lb]{\smash{{{\SetFigFont{11}{13.2}{\rmdefault}{\mddefault}{\itdefault}B}}}}}
\put(4859,2159){\makebox(0,0)[lb]{\smash{{{\SetFigFont{11}{13.2}{\rmdefault}{\mddefault}{\itdefault}B}}}}}
\put(6096,696){\makebox(0,0)[lb]{\smash{{{\SetFigFont{20}{24.0}{\rmdefault}{\mddefault}{\itdefault}C}}}}}
\put(247,6847){\makebox(0,0)[lb]{\smash{{{\SetFigFont{7}{8.4}{\rmdefault}{\mddefault}{\itdefault}A}}}}}
\put(8422,5797){\makebox(0,0)[lb]{\smash{{{\SetFigFont{7}{8.4}{\rmdefault}{\mddefault}{\itdefault}l}}}}}
\put(8497,5872){\makebox(0,0)[lb]{\smash{{{\SetFigFont{5}{6.0}{\rmdefault}{\mddefault}{\updefault}1}}}}}
\put(8497,5722){\makebox(0,0)[lb]{\smash{{{\SetFigFont{5}{6.0}{\rmdefault}{\mddefault}{\updefault}31}}}}}
\put(8947,6997){\makebox(0,0)[lb]{\smash{{{\SetFigFont{7}{8.4}{\rmdefault}{\mddefault}{\itdefault}l}}}}}
\put(9022,7072){\makebox(0,0)[lb]{\smash{{{\SetFigFont{5}{6.0}{\rmdefault}{\mddefault}{\updefault}1}}}}}
\put(9022,6922){\makebox(0,0)[lb]{\smash{{{\SetFigFont{5}{6.0}{\rmdefault}{\mddefault}{\updefault}12}}}}}
\put(9172,6697){\makebox(0,0)[lb]{\smash{{{\SetFigFont{7}{8.4}{\rmdefault}{\mddefault}{\itdefault}l}}}}}
\put(9247,6772){\makebox(0,0)[lb]{\smash{{{\SetFigFont{5}{6.0}{\rmdefault}{\mddefault}{\updefault}1}}}}}
\put(9547,6997){\makebox(0,0)[lb]{\smash{{{\SetFigFont{7}{8.4}{\rmdefault}{\mddefault}{\itdefault}l}}}}}
\put(9622,7072){\makebox(0,0)[lb]{\smash{{{\SetFigFont{5}{6.0}{\rmdefault}{\mddefault}{\updefault}1}}}}}
\put(9622,6922){\makebox(0,0)[lb]{\smash{{{\SetFigFont{5}{6.0}{\rmdefault}{\mddefault}{\updefault}13}}}}}
\put(9772,6697){\makebox(0,0)[lb]{\smash{{{\SetFigFont{7}{8.4}{\rmdefault}{\mddefault}{\itdefault}l}}}}}
\put(9847,6772){\makebox(0,0)[lb]{\smash{{{\SetFigFont{5}{6.0}{\rmdefault}{\mddefault}{\updefault}1}}}}}
\put(9847,6622){\makebox(0,0)[lb]{\smash{{{\SetFigFont{5}{6.0}{\rmdefault}{\mddefault}{\updefault}13}}}}}
\put(9847,6097){\makebox(0,0)[lb]{\smash{{{\SetFigFont{7}{8.4}{\rmdefault}{\mddefault}{\itdefault}l}}}}}
\put(9922,6172){\makebox(0,0)[lb]{\smash{{{\SetFigFont{5}{6.0}{\rmdefault}{\mddefault}{\updefault}1}}}}}
\put(9622,6397){\makebox(0,0)[lb]{\smash{{{\SetFigFont{7}{8.4}{\rmdefault}{\mddefault}{\itdefault}l}}}}}
\put(9697,6472){\makebox(0,0)[lb]{\smash{{{\SetFigFont{5}{6.0}{\rmdefault}{\mddefault}{\updefault}1}}}}}
\put(9697,6322){\makebox(0,0)[lb]{\smash{{{\SetFigFont{5}{6.0}{\rmdefault}{\mddefault}{\updefault}23}}}}}
\put(9022,5797){\makebox(0,0)[lb]{\smash{{{\SetFigFont{7}{8.4}{\rmdefault}{\mddefault}{\itdefault}l}}}}}
\put(9097,5872){\makebox(0,0)[lb]{\smash{{{\SetFigFont{5}{6.0}{\rmdefault}{\mddefault}{\updefault}1}}}}}
\put(9097,5722){\makebox(0,0)[lb]{\smash{{{\SetFigFont{5}{6.0}{\rmdefault}{\mddefault}{\updefault}32}}}}}
\put(9247,6097){\makebox(0,0)[lb]{\smash{{{\SetFigFont{7}{8.4}{\rmdefault}{\mddefault}{\itdefault}l}}}}}
\put(9322,6172){\makebox(0,0)[lb]{\smash{{{\SetFigFont{5}{6.0}{\rmdefault}{\mddefault}{\updefault}1}}}}}
\put(9322,6022){\makebox(0,0)[lb]{\smash{{{\SetFigFont{5}{6.0}{\rmdefault}{\mddefault}{\updefault}22}}}}}
\put(9022,6397){\makebox(0,0)[lb]{\smash{{{\SetFigFont{7}{8.4}{\rmdefault}{\mddefault}{\itdefault}l}}}}}
\put(9097,6472){\makebox(0,0)[lb]{\smash{{{\SetFigFont{5}{6.0}{\rmdefault}{\mddefault}{\updefault}1}}}}}
\put(9097,6322){\makebox(0,0)[lb]{\smash{{{\SetFigFont{5}{6.0}{\rmdefault}{\mddefault}{\updefault}22}}}}}
\put(8422,6397){\makebox(0,0)[lb]{\smash{{{\SetFigFont{7}{8.4}{\rmdefault}{\mddefault}{\itdefault}l}}}}}
\put(8497,6472){\makebox(0,0)[lb]{\smash{{{\SetFigFont{5}{6.0}{\rmdefault}{\mddefault}{\updefault}1}}}}}
\put(8497,6322){\makebox(0,0)[lb]{\smash{{{\SetFigFont{5}{6.0}{\rmdefault}{\mddefault}{\updefault}21}}}}}
\put(8647,6097){\makebox(0,0)[lb]{\smash{{{\SetFigFont{7}{8.4}{\rmdefault}{\mddefault}{\itdefault}l}}}}}
\put(8722,6172){\makebox(0,0)[lb]{\smash{{{\SetFigFont{5}{6.0}{\rmdefault}{\mddefault}{\updefault}1}}}}}
\put(8722,6022){\makebox(0,0)[lb]{\smash{{{\SetFigFont{5}{6.0}{\rmdefault}{\mddefault}{\updefault}21}}}}}
\put(8647,5497){\makebox(0,0)[lb]{\smash{{{\SetFigFont{7}{8.4}{\rmdefault}{\mddefault}{\itdefault}l}}}}}
\put(8722,5572){\makebox(0,0)[lb]{\smash{{{\SetFigFont{5}{6.0}{\rmdefault}{\mddefault}{\updefault}1}}}}}
\put(8722,5422){\makebox(0,0)[lb]{\smash{{{\SetFigFont{5}{6.0}{\rmdefault}{\mddefault}{\updefault}31}}}}}
\put(9247,5497){\makebox(0,0)[lb]{\smash{{{\SetFigFont{7}{8.4}{\rmdefault}{\mddefault}{\itdefault}32}}}}}
\put(9322,5572){\makebox(0,0)[lb]{\smash{{{\SetFigFont{5}{6.0}{\rmdefault}{\mddefault}{\updefault}1}}}}}
\put(9322,5422){\makebox(0,0)[lb]{\smash{{{\SetFigFont{5}{6.0}{\rmdefault}{\mddefault}{\updefault}11}}}}}
\put(9847,5497){\makebox(0,0)[lb]{\smash{{{\SetFigFont{7}{8.4}{\rmdefault}{\mddefault}{\itdefault}l}}}}}
\put(9922,5422){\makebox(0,0)[lb]{\smash{{{\SetFigFont{5}{6.0}{\rmdefault}{\mddefault}{\updefault}33}}}}}
\put(9622,5797){\makebox(0,0)[lb]{\smash{{{\SetFigFont{7}{8.4}{\rmdefault}{\mddefault}{\itdefault}l}}}}}
\put(9697,5872){\makebox(0,0)[lb]{\smash{{{\SetFigFont{5}{6.0}{\rmdefault}{\mddefault}{\updefault}1}}}}}
\put(9697,5722){\makebox(0,0)[lb]{\smash{{{\SetFigFont{5}{6.0}{\rmdefault}{\mddefault}{\updefault}33}}}}}
\put(9247,6622){\makebox(0,0)[lb]{\smash{{{\SetFigFont{5}{6.0}{\rmdefault}{\mddefault}{\updefault}12}}}}}
\put(9922,6022){\makebox(0,0)[lb]{\smash{{{\SetFigFont{5}{6.0}{\rmdefault}{\mddefault}{\updefault}23}}}}}
\put(9922,5572){\makebox(0,0)[lb]{\smash{{{\SetFigFont{5}{6.0}{\rmdefault}{\mddefault}{\updefault}1}}}}}
\end{picture}
}
\caption{{}$\hat{A}\subset B(H)$ and $L$. Here the dotted squares are matrix units of $B(H)$; 
 large solid-line squares are matrix units in $B(K_{j})$; $A\in B(H_{1})$, $B\in B(H_{2})$, 
   $C\in B(H_{3})$, and $d(1)=3$, $d(2)=4$ and $d(3)=1$.}\label{figure:AandL}\end{figure}

Let
\[
L_{p}=\sum _{k}\sum _{i,j=1}^{d(k)}\sqrt{\lambda _{i}(k)}l_{ij}^{k}(p)\otimes (1_{B(H_{k})}\otimes e_{ij}(k))\in P\]
(see Figure \ref{figure:AandL}). Then for all \( T\in B(H) \), we get
\begin{eqnarray*}
L_{p}^{*}TL_{q} & = & \delta _{pq}\left( \sum _{k}\sum _{i,j=1}^{d(k)}\sqrt{\lambda _{i}(k)}l_{ij}^{k}(p)\otimes 
(1_{B(H_{k})}\otimes e_{ij}(k))\right) ^{*}T\\
 &  & \cdot \sum _{k'}\sum _{i',j'=1}^{d(k')}\sqrt{\lambda _{i'}(k')}l_{i'j'}^{k'}(p)\otimes (1_{B(H_{k'})}\otimes e_{i'j'}(k))
\\
 & = & \delta _{pq}\sum _{k}\sum _{i,j}^{d(k)}\lambda _{i}(k)(1_{B(H_{k})}\otimes e_{ji}(k))T(1_{B(H_{k})}\otimes e_{ij}(k))\\
 & = & \delta _{pq}E(T).
\end{eqnarray*}
Hence by \cite[Theorem 2.3]{shlyakht:amalg}, \( \{L_{p}:p=1,2,\ldots \} \)
are \( * \)-free with respect to \( \psi\otimes E :P\to B(H) \) from \( B(H) \)
with amalgamation over \( \hat{A} \). Moreover, since \( L_{p}^{*}aL_{q}=\delta _{pq}a \)
for all \( a\in \hat{A} \), we get that \( W^{*}(\hat{A},L_{p}+L_{p}^{*}:p\geq 1)\cong W^{*}(L_{p}+L_{p}^{*}:p\geq 1)\otimes 
\hat{A}\cong L(\mathbb {F}_{\infty })\otimes \hat{A} \)
(cf. \cite[Examples 3.2 and 3.3(b)]{shlyakht:semicirc}). We conclude that \( M=W^{*}(B(H),L_{p}+L_{p}^{*})\cong (L(\mathbb 
{F}_{\infty })\otimes \hat{A},\tau \otimes \id )*_{\hat{A}}(B(H),E) \). 

Now, fix a minimal projection \( p\in B(H_{1})\otimes e_{11}(1) \). Then the
compression \( pMp \) is generated as a von Neumann algebra by the entries
of \( \{L_{p}+L_{p}^{*}:p\geq 1\} \) viewed as matrices. These entries have
the form \( \sqrt{\lambda _{i}(k)}l_{ij}^{k}(p)+\sqrt{\lambda _{j}(k)}(l_{ji}^{k}(p))^{*} \),
\( k=1,2,\ldots  \), \( 1\leq i\leq j\leq d(k) \), \( p\geq 1 \). The algebra
generated by such elements is isomorphic to a free Araki-Woods factor (see \cite[Section 5]{shlyakht:quasifree:big});
the classifying \( Sd \) invariant for this type of factors is the multiplicative
subgroup of \( \mathbb {R} \) generated by the ratios \( \lambda _{i}(k)/\lambda _{j}(k) \),
\( k=1,2,\ldots  \), \( 1\leq i,j\leq d(k) \).

We now turn to the particular case of \( \hat{A} \) arising from a quantum
group action, where one can describe the inclusion \( B(L^{2}(A))\supseteq \hat{A} \)
and the conditional expectation \( E \) explicitly. We omit the details of
the computation, but give only the consequence for the reader's convenience.

Let \( w^{(\ell )}_{ij} \), \( i,j\in I_{\ell }=\{-\ell ,\ell +1,\dots ,\ell \} \),
be the matrix elements of the spin \( \ell  \) (\( \in (1/2){\mathbb N}_{0} \))
representation as given in \cite{masudaANDCO:representationsSUQ2andQJacobiPolynomials}.
The Peter-Weyl type theorem says that the set
\[
\left\{ \left( \frac{1-q^{2(2\ell +1)}}{(1-q^{2})q^{2(\ell +i)}}\right) ^{1/2}\Lambda _{h}(w_{ij}^{\ell })\right\} _{i,j\in 
I_{\ell },\ell \in (1/2){\mathbb N}_{0}}\]
forms a complete orthonormal basis of \( L^{2}(A) \). Let us denote the matrix
units with respect to the basis by \( e_{(i_{1},j_{1},\ell _{1})(i_{2},j_{2},\ell _{2})} \)'s.
The set of elements of the form:
\[
\hat{E}^{(\ell )}_{i_{1}i_{2}}=\sum _{j\in I_{\ell }}e_{(i_{1},j,\ell )(i_{2},j,\ell )}\]
forms a complete system of matrix units of \( \hat{A} \), and the conditional
expectation \( E \) is computed as follows:
\[
E\left( e_{(i_{1},j_{1},\ell _{1})(i_{2},j_{2},\ell _{2})}\right) =\delta _{(i_{1},\ell _{1})(i_{2},\ell _{2})}\frac{1}{1+q^{2}
+\cdots +(q^{2})^{2\ell _{1}}}q^{2(\ell _{1}-i_{1})}\hat{E}^{(\ell _{1})}_{j_{1}j_{2}}.\]
In other words, in the notation of the first part of the appendix, we have \( \dim H_{\ell }=\dim K_{\ell }=d(\ell )=2\ell +1 
\),
and the eigenvalues of \( \theta _{\ell } \) are given by
\[
\lambda _{i}(\ell )=\frac{1}{1+q^{2}+\cdots +(q^{2})^{2\ell }}q^{2(\ell -i)}.\]
Since the \( Sd \) invariant of the associated type III factor is the multiplicative
subgroup of positive reals generated by the ratios \( \lambda _{i}(\ell )/\lambda _{j}(\ell ) \),
we see that in this case \( Sd=q^{2\mathbb {Z}} \), i.e., the factor is of type
III\( _{q^{2}} \).

\section*{Appendix II.}

For reader's convenience, we give here a direct argument showing that 
the conditional expectation $E_{B(H)} : \Phi(B(H), \eta) \rightarrow B(H)$ 
in Theorem 2.2 is faithful. We use the notations given in the beginning 
of \S2. For more on this, interested readers are refered to 
\cite{shlyakht:semicirc}. 

\begin{lem*}
\label{lemma:KMScondition}{\rm (cf. \cite[Lemma 4.2]{shlyakht:amalg})} 
Suppose that 
\( \eta_{ij} = \delta_{ij} E \), where \( E : M \rightarrow B \) is 
a faithful normal conditional expectation. Fix a faithful normal state
\( \varphi \) on \( B \). Then the state 
\( \theta := \varphi\circ E\circ E_M \) on \( \Phi(M, \eta) \) 
satisfies 
\begin{eqnarray*}
\theta(f_0 X_{i_1} f_1 \cdots X_{i_{n-1}} f_{n-1} X_{i_n} f_n) 
& = & \theta(\sigma_{i}^{\varphi\circ E}(f_n) f _0 X_{i_1} f_1 \cdots X_{i_{n-1}} f_{n-1} X_{i_n}) \\
& = & \theta(X_{i_n} \sigma_{i}^{\varphi\circ E}(f_n) f _0 X_{i_1} f_1 \cdots  X_{i_{n-1}} f_{n-1}) 
\end{eqnarray*}
for all analytic elements \( f_0, f_1,\dots, f_n \in M \) with respect to 
the modular action \( \sigma^{\varphi\circ E} \). 
\end{lem*} 
\begin{proof} The proof is a straightforward modification of 
\cite[Lemma 4.2]{shlyakht:amalg}. 
Note that \( \varphi\circ E(x E(y)) 
= \varphi\circ E(E(x) y) \) (\( x, y \in M \)),
so we use the following identities (instead of the trace property of 
\(\varphi\circ E \) 
used in Lemma 4.2):  \\
(1) \( \sigma^{\varphi\circ E}_z\circ E(x) = 
E\circ\sigma_z^{\varphi\circ E}(x) \) for analytic \( x \in M \) and 
\( z \in {\mathbb C} \); \\ 
(2) for analytic \( x, y \in M \)  
\[
\theta(xy) = \varphi\circ E(xy) = 
\varphi\circ E(\sigma_{i}^{\varphi\circ E}(y)x) = 
\theta(\sigma_{i}^{\varphi\circ E}(y)x).  
\]
\end{proof}

Since \( \sigma_t^{\varphi\circ E}\circ E = E\circ\sigma_t^{\varphi\circ
E} \) and \( \varphi\circ E \circ\sigma_t^{\varphi\circ E} = \varphi\circ
E \) (\( t \in {\mathbb R} \)),  there exists a 1-parameter
group of automorphisms of \(\Phi(M,\eta)\) 
\( \sigma_t \) (\( t \in {\mathbb R} \)),  for which
\[
\sigma_t|_M = \sigma_t^{\varphi\circ E}, \quad \sigma_t(X_i) = X_i 
\]
(note that \( \theta\circ\sigma_t = \theta \)).
Hence the above lemma implies the following: 

\begin{prop*}
\label{prop:A-valuedmodular} With the assuputions and notations in 
the Lemma, 
the action \( \sigma_t \) {\rm (}\( t \in {\mathbb R} \){\rm )} satisfies 
the modular condition for 
\( \theta \), i.e., \( \sigma_t = \sigma_t^{\theta} \), 
and hence \( \theta \) is faithful. In particular, the canonical
conditional expectation \( E_M : \Phi(M, \eta) \rightarrow M \) is faithful. 
\end{prop*} 

\bibliographystyle{amsplain}

\providecommand{\bysame}{\leavevmode\hbox to3em{\hrulefill}\thinspace}

\end{document}